\documentstyle[12pt]{article}
\newtheorem{theo}{Theorem}

\newtheorem{prop}[theo]{Proposition}
\newtheorem{lemm}[theo]{Lemma}
\def\al{\alpha}
\def\be{\beta}
\def\ga{\gamma}
\def\de{\delta}

\def\la{\lambda}

\def\De{\Delta}
\def\beq{\begin{equation}}
\def\eeq{\end{equation}}
\def\bea{\begin{eqnarray}}
\def\eea{\end{eqnarray}}
\def\beas{\begin{eqnarray*}}
\def\eeas{\end{eqnarray*}}
\def\nn{\nonumber}
\def\ds{\displaystyle}
\def\U{{\cal U}_h(sl(2))}
\renewcommand{\theequation}{\arabic{section}.\arabic{equation}}
\begin{document}
\begin{center}
{\Large \bf Representations and Clebsch-Gordan coefficients}\\[2mm]
{\Large \bf for the Jordanian quantum algebra ${\cal
U}_h(sl(2))$}\\[2cm] 
J.~Van der Jeugt\footnote{Senior Research Associate  of the Fund
for Scientific Research -- Flanders (Belgium)} \\[8mm]
Vakgroep Toegepaste Wiskunde en Informatica,
Universiteit Gent,\\
Krijgslaan 281--S9, B9000 Gent, Belgium
\end{center}

\vskip 3cm
\begin{abstract}
Representation theory for the Jordanian quantum algebra $\U$ is
developed. Closed form expressions are given for the action of the
generators of $\U$ on the basis vectors of finite dimensional
irreducible representations. It is shown how representation theory of
$\U$ has a close connection to combinatorial identities involving
summation formulas.
A general formula is obtained for the
Clebsch-Gordan coefficients ${\cal C}^{j_1,j_2,j}_{n_1,n_2,m}(h)$ of
$\U$. These Clebsch-Gordan coefficients are shown to coincide with
those of $su(2)$ for $n_1+n_2 \leq m$, but for $n_1+n_2 > m$ they are
in general a nonzero monomial in $h^{n_1+n_2-m}$. 
\end{abstract}

\newpage
\section{Introduction}
\setcounter{equation}{0}

Viewing a quantum group as a quantum automorphism group acting on a
noncommuting space~\cite{Manin,Schupp,Dubois}, one often requires the
extra condition 
of the existence of a central determinant in the quantum matrix group.
For two by two matrices, this condition restricts the quantum groups
essentially to only two classes~\cite{Kuperschmidt},
namely the standard $SL_q(2)$ quantum group and the Jordanian $SL_h(2)$
quantum group. The quantum group $SL_h(2)$ was introduced
in~\cite{Zakrzewski}, and the corresponding quantum algebra (or
quantised universal enveloping algebra) $\U$ was given in~\cite{Ohn}. A
unversal $R$-matrix for $\U$ was constructed in~\cite{Ballesteros}.

The main object of this paper is to develop representation theory of
the Jordanian quantum algebra $\U$, and in particular construct
Clebsch-Gordan coefficients. The finite dimensional highest weight
representations of $\U$ were given in~\cite{Dobrev}, first by a direct
construction, and then by factorising the corresponding Verma module.
In~\cite{Dobrev} the action of the $\U$ generators on a finite dimensional
representation was not given explicitly. An important construction was
developed by Abdesselam {\em et al}~\cite{Abdesselam}~: they gave a
nonlinear relation between the generators of $\U$ and the classical
generators of $sl(2)$. As a consequence of this relation, they obtained
expressions for the action of the $\U$ generators $H$, $X$ and $Y$ (see
the following section for their definition) on basis vectors of the
finite dimensional irreducible representations. These expressions are
in closed form, except for the action of the generator $Y$. Using this
nonlinear map, Aizawa~\cite{Aizawa} constructed finite and infinite
dimensional representations of $\U$, and considered the tensor product
of two representations. Moreover, he gives some examples of
Clebsch-Gordan coefficients. 

Our present work is motivated by the fact that $sl(2)$ or $su(2)$
representations appear in many physical theories, and often their
Clebsch-Gordan coefficients are fundamental in these theories. Since
representation theory of $\U$ is so closely related to that of $su(2)$,
and could be used as the algebraic structure underlying deformations of
such physical models, an important aspect to study are the Clebsch-Gordan
coefficients of $\U$. In the present paper it is shown how an explicit
formula for $\U$ Clebsch-Gordan coefficients can be obtained.

In Section~2, we give the defining relations for $\U$, and the
nonlinear relation between $sl(2)$ generators and $\U$ generators. In
Section~3, closed forms for the action of the three generators $H$, $X$
and $Y$ of $\U$ acting on basis vectors of finite dimensional
irreducible representations are determined. For $H$ and $X$, these
expressions correspond to those of Ref.~\cite{Abdesselam}; the
determination of the explicit action of $Y$ is new and is found using a
number of combinatorial identities (Lemmas~\ref{lemmYcoef}
and~\ref{lemmzeilb}). In Section~4 the tensor product of two
representations is considered. In this tensor product we show the
existence of an auxiliary basis which behaves like the uncoupled basis
vectors in the tensor product of two $su(2)$ representations. Using
this auxiliary basis, the $\U$ Clebsch-Gordan coefficients are easily
determined in Section~5, and some expamles and properties are discussed
in Section~6. 

A curious aspect of the results in this paper (and more generally of
$\U$ representation theory) is that they are closely related to
nontrivial combinatorial identities (see Lemmas~\ref{lemmYcoef},
\ref{lemmzeilb} and~\ref{lemmgosper}). The identities needed here have
on one side a (definite or indefinite) sum of hypergeometric terms, and
a closed form expression on the other side. To find closed form
expressions for such summations is a problem that can be solved
completely algorithmically~\cite{AB}~: for indefinite summations this
can be done by means of Gosper's algorithm; for definite summations
this is done by means of Zeilberger's algorithm. For both algorithms,
programs are available in Maple or Mathematica. In the Appendix, we
comment on the proofs of these combinatorial identities.

\section{Definition and relation to $sl(2)$}
\setcounter{equation}{0}

The Jordanian quantum algebra $\U$ is an associative algebra with unity
and generated by $X$, $Y$ and $H$, subject to the relations
\bea
&&[H,X]=2 {\sinh hX \over h}, \qquad 
[H,Y]= -Y(\cosh hX) - (\cosh hX)Y, \nn \\
&&[X,Y]=H;
\eea
herein, $h$ is the deformation parameter. We shall assume that $|h|<1$.
In the limit $h\rightarrow 0$, $\U$ reduces to the universal enveloping
algebra of $sl(2)$. The Hopf algebra structure of $\U$ is given
in~\cite{Ohn}; here, we are only interested in the comultiplication, which
reads
\bea
&&\De(X)=X\otimes 1 + 1\otimes X,\nn \\
&&\De(Y)=Y\otimes e^{hX}+ e^{-hX}\otimes Y, \label{Del}\\
&&\De(H)=H\otimes e^{hX}+ e^{-hX}\otimes H. \nn
\eea

The irreducible finite-dimensional highest weight representations of
$\U$ can be obtained by using the invertible map from $sl(2)$ to $\U$,
given in~\cite{Abdesselam}. With the following definition
\bea
&& Z_+ = {2 \over h} \tanh {hX \over 2},\nn\\
&& Z_- = ( \cosh {hX \over 2}) Y (\cosh {hX\over 2}), 
\label{sl2toU}
\eea
it follows that the elements $\{H, Z_+, Z_-\}$ satisfy the commutation
relations of a classical $sl(2)$ basis~:
\beq
[H,Z_\pm]=\pm 2 Z_{\pm}, \qquad [Z_+,Z_-]=H.
\label{defrelsl2}
\eeq
The relations (\ref{sl2toU}) can be inverted, see the following
section, and thus with every finite-dimensional irreducible $sl(2)$
representation there corresponds a finite-dimensional irreducible
representation of $\U$. These $sl(2)$ representations are labeled by a
number $j$, with $2j$ a non-negative integer; the representation space
is $V^{(j)}$ with basis $e^j_m$, where $m= -j, -j+1,\ldots, j$. The
action of $sl(2)$ on this basis is given by
\bea
&& H e^j_m = 2m\; e^j_m, \nn\\
&& Z_\pm e^j_m = \sqrt{ (j\mp m)(j\pm m+1) }\; e^j_{m\pm 1}.
\label{sl2eaction}
\eea
For most of the computations in this paper, it is easier to work with
another basis of $V^{(j)}$ related to the above basis by
\beq
v^j_m = \al_{j,m} e^j_m,\;\;\hbox{with } \al_{j,m} = \sqrt{
(j+m)!/(j-m)!}\; .
\label{ve}
\eeq
The matrix elements of the $sl(2)$ generators is then given by
\bea
&& H v^j_m = 2m\; v^j_m, \nn\\
&& Z_+ v^j_m = v^j_{m+1},  \label{sl2vaction} \\
&& Z_- v^j_m = (j+ m)(j- m+1)\; v^j_{m-1}, \nn
\eea
where $v^j_{j+1}=0$. Clearly, if for an operator the matrix
elements in the $v$-basis have been determined, the matrix elements in
the $e$-basis follow immediately using~(\ref{ve}).

\section{Representations of $\U$}
\setcounter{equation}{0}

In this section we wish to give explicit expressions for the matrix
elements of $H$, $X$ and $Y$ in the $v$-basis. For $H$, this is
trivial, see~(\ref{sl2vaction}). For $X$, one first determines the
action of $e^{hX}$. From relation~(\ref{sl2toU}) one finds that 
\beq
e^{hX} = (1+{h\over 2} Z_+)   (1-{h\over 2} Z_+)^{-1}.
\label{XitoZ}
\eeq
Then the action of $Z_+$ in the $v$-basis implies
\beq
e^{hX} v^j_m = v^j_m + 2 \sum_{k=1}^{j-m} ({h\over 2})^k v^j_{m+k}.
\label{actionehX}
\eeq
Thus in this representation one can write $e^{hX}=1+N_1$, with $N_1$ a
nilpotent matrix. Then $hX=\log(1+N_1)=N_1-N_1^2/2 +N_1^3/3 -\cdots$,
and one obtains the following action of $X$ in the representation space
$V^{(j)}$~: 
\beq
X v^j_m = \sum_{k=0}^{\lfloor (j-m-1)/2 \rfloor} {(h/2)^{2k} \over
2k+1} \; v^j_{m+1+2k}.
\label{Xaction}
\eeq
Up to a scaling of the basis vectors, (\ref{Xaction}) coincides with
\cite[eq.~(23)]{Abdesselam}. 
The action of $Y$ is more difficult to determine explicitly; in
\cite[eq.~(35)]{Abdesselam} an expression is given but the matrix
elements still involve a single sum. Here we shall show that it can
actually be given in closed form. Let us use the relation 
\beq
 Y = ( \cosh {hX \over 2})^{-1} Z_- (\cosh {hX\over 2})^{-1},
\label{YitoZ}
\eeq
and first determine the matrix form of $(\cosh {hX\over 2})^{-1}$ in
the $v$-basis. In this basis of $V^{(j)}$, one can write
\beq
(\cosh {hX\over 2})^2 = {1\over 2} (1 + \cosh hX)
= {1\over 2} (1 + {e^{hX}+e^{-hX} \over 2}) = 1+N_2,
\eeq
where, by (\ref{actionehX}), $N_2$ is again a nilpotent matrix whose
matrix elements follow from those of $e^{hX}$~:
\beq
N_2 v^j_m = \sum_{k=1}^{\lfloor(j-m)/2 \rfloor} (h/2)^{2k}
v^j_{m+2k}. 
\eeq
Then, in this representation
\bea
(\cosh {hX\over 2})^{-1} &=& (1+N_2)^{-1/2} \nn\\
&=& \sum_{n=0}^\infty (-1)^n (1/2)_n {N_2^n \over n!},
\label{expN2}
\eea
where $(a)_n$ is the notation for the Pochhammer symbol~:
\beq
(a)_n = \left\{ 
 \begin{array}{lll} 
 a(a+1)\cdots (a+n-1) & \hbox{if} & n=1,2,\ldots; \\
 1 & \hbox{if} & n=0.
 \end{array} \right.
\label{poch}
\eeq
Usin the explicit powers of $N_2$ in the action of (\ref{expN2}) on
$v^j_m$, the contributions to the
coefficient of $v^j_{m+2k}$ ($k>0$) in $(1+N_2)^{-1/2} v^j_m$ are 
\beq
\left( h\over 2 \right)^{2k} \sum_{n=1}^{k} (-1)^n {(1/2)_n \over n!}
\left({k-1 \atop n-1}\right).
\eeq

Next, we use
\begin{lemm}
For $k>0$ integer,
$$
\sum_{n=1}^{k} (-1)^n {(1/2)_n \over n!} \left({k-1 \atop n-1}\right) =
- {1\over 2^{2k-1}} {(2k-2)! \over k!(k-1)!}.
$$
\label{lemmYcoef}
\end{lemm}

As a consequence, one finds the explicit action of $(\cosh {hX\over
2})^{-1}$ in the $v$-basis~:
\beq
(\cosh {hX\over 2})^{-1} v^j_m =
v^j_m -2 \sum_{k=1}^{\lfloor(j-m)/2
\rfloor} t_k \left(h\over 4\right)^{2k} v^j_{m+2k},
\label{cosh-1}
\eeq
with $t_k={(2k-2)!\over k!(k-1)!}$. Using (\ref{YitoZ}), (\ref{cosh-1})
and (\ref{sl2vaction}), one determines the action of $Y$,
\bea
Yv^j_m &=& (\cosh {hX\over 2})^{-1} Z_- \left(v^j_m -2 \sum_{k\geq 1}
t_k (h/4)^{2k} v^j_{m+2k} \right) \nn\\
&=& (\cosh {hX\over 2})^{-1} \Bigl((j+m)(j-m+1) v^j_{m-1} - \nn\\ 
&&  2 \sum_{k\geq 1}
t_k (h/4)^{2k} (j+m+2k)(j-m-2k+1) v^j_{m+2k-1} \Bigr)\nn\\
&=& (j+m)(j-m+1) \left(v^j_{m-1} -2 \sum_{l\geq 1}
t_l (h/4)^{2l} v^j_{m+2l-1} \right) - \nn\\
&& 2 \sum_{k\geq 1} t_k (h/4)^{2k} (j+m+2k)(j-m-2k+1)\nn\\
&&\qquad\times \left(
v^j_{m+2k-1} -2 \sum_{l\geq 1} t_l (h/4)^{2l} v^j_{m+2k+2l-1} \right) .
\label{tempresY}
\eea
In this last expression, the coefficient of $v^j_{m-1}$ is
$(j+m)(j-m+1)$. The coefficient of $v^j_{m+1}$ is
\bea
&&-2(j+m)(j-m+1)t_1 h^2/16 - 2(j+m+2)(j-m-1)t_1 h^2/16  \nn\\
&&= h^2/4 -(j-m)(j+m+1) h^2/4.
\eea
Finally, we determine the coefficient of $v^j_{m+2s-1}$ ($s\geq 2$) in
(\ref{tempresY}); it reads
\bea
&&-2(j+m)(j-m+1)t_s (h/4)^{2s} -2(j+m+2s)(j-m-2s+1) t_s (h/4)^{2s} +
\nn\\ 
&& \qquad 4 (h/4)^{2s} \sum_{k=1}^{s-1} (j+m+2k)(j-m-2k+1) t_k t_{s-k}.
\label{coeffv}
\eea
In order to give a closed form for this coefficient, we
need to find an expression for $\sum_k t_k t_{s-k} k^n$ with $n=0$, 1
and 2. This is done with the help of the following lemma~:

\begin{lemm}
Let $s\geq 2$ be integer, then
$$
\sum_{k=1}^{s-1} {(2k-2)!\over k!(k-1)!}
{(2s-2k-2)!\over(s-k)!(s-k-1)!} k^n =
\left\{ 
 \begin{array}{ll}
 \ds{(2s-2)! \over s!(s-1)!} & \hbox{if } n=0,\\[4mm]
 \ds{(2s-2)! \over 2(s-1)!^2} & \hbox{if } n=1,\\[4mm]
 \ds{s(2s-2)! \over 2(s-1)!^2}-4^{s-2} & \hbox{if } n=2.
 \end{array}
\right.
$$
\label{lemmzeilb}
\end{lemm}

Using these three results, it follows from (\ref{coeffv}) that the
coefficient of $v^j_{m+2s-1}$ ($s\geq 2$) is $4^s
(h/4)^{2s}=(h/2)^{2s}$. Thus we finally obtain the explicit action of
$Y$ on the $v$-basis~:
\bea
Y v^j_m &=& (j+m)(j-m+1) v^j_{m-1} -(j-m)(j+m+1)\left(h\over 2\right)^2
v^j_{m+1} + \nn\\
&&\qquad \sum_{s=1}^{\lfloor(j-m+1)/2\rfloor} \left(h\over
2\right)^{2s} v^j_{m-1+2s}.
\label{Yaction}
\eea
Thus the first equation of (\ref{sl2vaction}) together with
equations~(\ref{Xaction}) and~(\ref{Yaction}) determine the action of
the generators $H$, $X$ and $Y$ of $\U$ on the finite dimensional
irreducible representations in closed form.

\section{Tensor product of two representations}
\setcounter{equation}{0}

The action of any of the generators $H$, $X$ or $Y$ on the tensor
product of two representations is determined by the comultiplication.
The comultiplication rule on $H$, $X$ and $Y$ induces a
comultiplication on $H$, $Z_\pm$. The purpose of this section is to
show that the tensor product $V^{(j_1)}\otimes V^{(j_2)}$ decomposes
into a direct sum of representations isomorphic to $V^{(j)}$, where
$j=j_1+j_2, j_1+j_2-1, \ldots, |j_1-j_2|$, and to determine the
Clebsch-Gordan coefficients in
\beq
e^{(j_1j_2)j}_m = \sum_{n_1,n_2} {\cal C}^{j_1,j_2,j}_{n_1,n_2,m}(h) \;
e^{j_1}_{n_1} \otimes e^{j_2}_{n_2},
\eeq
such that $e^{(j_1j_2)j}_m$ is a standard $e$-basis of $V^{(j)}$, i.e.
\bea
&& \De(H) e^{(j_1j_2)j}_m = 2m\; e^{(j_1j_2)j}_m, \nn\\
&& \De(Z_\pm) e^{(j_1j_2)j}_m = \sqrt{ (j\mp m)(j\pm m+1) }\;
e^{(j_1j_2)j}_{m\pm 1}. 
\eea

The first step towards this goal is to find expressions of $\De(H)$ and
$\De(Z_+)$ in terms of $H$ and $Z_+$. For $H$, the problem is easy, and
follows from (\ref{Del}) and (\ref{XitoZ}), see
also~\cite[eq.~(3.2)]{Aizawa}~: 
\bea
&&\De(H) = H\otimes e^{hX} + e^{-hX} \otimes H \nn\\
&&= H\otimes 1 + 1 \otimes H + 2H \otimes \sum_{n=1}^\infty
\left(hZ_+\over 2 \right)^n +  \sum_{n=1}^\infty
\left(- hZ_+\over 2 \right)^n \otimes 2H.
\label{DelH}
\eea

To determine $\De(Z_+)$, denote $t=e^{hX}$ and $z={h\over 2} Z_+$. Then
$z=(t-1)(t+1)^{-1}$ and $t=(1+z)(1-z)^{-1}$. Moreover,
\beq
\De(t)=t\otimes t = (1+2\sum_{k=1}^\infty z^k) \otimes
(1+2\sum_{k=1}^\infty z^k). 
\eeq
Thus $\De(z)=\sum_{k,l=0}^\infty \la_{k,l} \; z^k\otimes z^l$, and the
coefficients can be obtained from 
$$
( 1\otimes 1 - \De(z)) \De(t) = 1\otimes 1 + \De(z),
$$
which must hold since $(1-z)t=(1+z)$. This leads to
$$
\De(z) = (1\otimes z + z \otimes 1)(1\otimes 1 - z\otimes z +
z^2\otimes z^2 - \cdots),
$$
and we obtain~:
\beq
\De(Z_+) = (1\otimes Z_+ + Z_+ \otimes 1)\left(\sum_{n=0}^\infty (-h^2/4)^n
 \; Z_+^n \otimes Z_+^n \right). 
\label{DelZ}
\eeq

An explicit expression for $\De(Z_-)$ is much more complicated (see,
e.g.\ eq.~(5.3) of~\cite{Aizawa}), but we do not need it here.

In the tensor product space $V^{(j_1)}\otimes V^{(j_2)}$, we now
define an auxiliary basis that is expressed in terms of the
$v$-basis $v^{j_1}_{m_1} \otimes v^{j_2}_{m_2}$. Define the
coefficients (see~(\ref{poch}) for the notation of Pochhammer symbols) 
\beq
b^{m_1,m_2}_{k,l} = 
\left\{ 
 \begin{array}{ll}
\ds{(-2m_1-k)_l (-2m_2-l)_k \over k!l!} & \hbox{if }k\geq 0\hbox{ and }
l\geq 0 ; \\[2mm]
0 & \hbox{otherwise},
 \end{array} \right.
\label{defb}
\eeq
and 
\beq
a^{m_1,m_2}_{k,l} = (-1)^k (h/2)^{k+l} (b^{m_1,m_2}_{k,l} -
 b^{m_1,m_2}_{k-1,l-1} ).
\label{defa}
\eeq
The auxiliary vectors are defined as follows~:
\beq
w^{j_1,j_2}_{m_1,m_2} = \sum_{k=0}^{j_1-m_1} \sum_{l=0}^{j_2-m_2}
a^{m_1,m_2}_{k,l} \; v^{j_1}_{m_1+k} \otimes v^{j_2}_{m_2+l}.
\label{defw}
\eeq
Clearly, they also form a basis for $V^{(j_1)}\otimes V^{(j_2)}$ since
the relation between the vectors $w^{j_1,j_2}_{m_1,m_2}$ and the
vectors $v^{j_1}_{m_1} \otimes v^{j_2}_{m_2}$ is given by an upper
triangular matrix with 1's on the diagonal. This auxiliary basis has
been introduced because the action of $\De(H)$ and $\De(Z_\pm)$ on it
is simple. The idea of introducing such an auxiliary basis comes
{}from~\cite{Aizawa}~: however, the coefficients used
in~\cite[eq.~(3.9)]{Aizawa} are single sum expressions and no closed
forms. 

\begin{prop} The action of $\De(H)$ on the auxiliary basis vectors is
given by
$$
\De(H) w^{j_1,j_2}_{m_1,m_2} = 2(m_1+m_2) \; w^{j_1,j_2}_{m_1,m_2}.
$$
\label{propH}
\end{prop}

To prove this, use (\ref{DelH}), (\ref{defw}), and the explicit action
of $H$ and $Z_+$ on the $v$-basis~:
\bea
&&\De(H) w^{j_1,j_2}_{m_1,m_2} = \De(H) \left(\sum_{k,l\geq 0} 
a^{m_1,m_2}_{k,l} \; v^{j_1}_{m_1+k} \otimes v^{j_2}_{m_2+l} \right) \nn\\
&&=\sum_{k,l\geq 0} a^{m_1,m_2}_{k,l} \Bigl( 2(m_1+k+m_2+l) 
v^{j_1}_{m_1+k} \otimes v^{j_2}_{m_2+l} \; +  \nn\\
&& 4(m_1+k) v^{j_1}_{m_1+k} \otimes \sum_{n\geq 1} (h/2)^n
v^{j_2}_{m_2+l+n} + 4(m_2+l) \sum_{n\geq 1} (-h/2)^n v^{j_1}_{m_1+k+n}
\otimes v^{j_2}_{m_2+l} \Bigr) \nn\\
&&= 2(m_1+m_2) w^{j_1,j_2}_{m_1,m_2}  + 2 \sum_{k,l\geq 0} 
v^{j_1}_{m_1+k} \otimes v^{j_2}_{m_2+l} \Bigl( (k+l) a^{m_1,m_2}_{k,l}
+  \nn\\
&&\quad 2(m_1+k) \sum_{n=1}^l (h/2)^n a^{m_1,m_2}_{k,l-n} +
2(m_2+l) \sum_{n=1}^l (-h/2)^n a^{m_1,m_2}_{k-n,l} \Bigr).
\label{DelHw}
\eea
Thus the proposition is proved provided we can show that the
coefficient of $v^{j_1}_{m_1+k} \otimes v^{j_2}_{m_2+l}$ in the last
summation is zero. Consider first
\beq
\sum_{n=1}^l (h/2)^n a^{m_1,m_2}_{k,l-n} = (-1)^k (h/2)^{k+l}
\sum_{n=1}^l (b^{m_1,m_2}_{k,l-n} - b^{m_1,m_2}_{k-1,l-n-1});
\eeq
we shall show that this last sum can be performed and yields
\beq
\sum_{n=1}^l (b^{m_1,m_2}_{k,l-n} - b^{m_1,m_2}_{k-1,l-n-1}) =
{ 2m_1+k-l+1 \over 2m_1+k } b^{m_1,m_2}_{k,l-1}.
\label{resb}
\eeq
Using the symmetry $b^{m_1,m_2}_{k,l}=b^{m_2,m_1}_{l,k}$, one can use
the same result (\ref{resb}) to find
\beq
\sum_{n=1}^l (-h/2)^n a^{m_1,m_2}_{k-n,l} = (-1)^k (h/2)^{k+l}
\Bigl( { 2m_2+l-k+1 \over 2m_2+l }\Bigr) b^{m_1,m_2}_{k-1,l}.
\eeq
Then the coefficient of $v^{j_1}_{m_1+k} \otimes v^{j_2}_{m_2+l}$ in
the summation-part of (\ref{DelHw}) is, up to a factor $(-1)^k
(h/2)^{k+l}$, equal to 
\bea
&&(k+l)(b^{m_1,m_2}_{k,l} - b^{m_1,m_2}_{k-1,l-1}) + 2(m_1+k) 
\Bigl({ 2m_1+k-l+1 \over 2m_1+k }\Bigr) b^{m_1,m_2}_{k,l-1} + \nn\\
&&\qquad 2(m_2+l) \Bigl({2m_2+l-k+1 \over 2m_2+l }\Bigr)
b^{m_1,m_2}_{k-1,l}, 
\eea
and using the definition of the coefficients $b^{m_1,m_2}_{k,l}$ this
is trivially shown to be zero. What remains to be proved is
(\ref{resb}). This follows from the following lemma~:

\begin{lemm}
Let $\al,\be,\ga$ be arbitrary parameters, $l\geq 0$ integer, and
$$
g(l)= {(\al)_l(\be)_l \over (l+1)! (\ga)_l } \bigl(
(\ga-\al-\be)l+\ga-1-\al\be \bigr). 
$$
Then for $n\geq 0$ integer,
$$
\sum_{l=0}^n g(l) = \ga-1-{ (\al)_{n+1}(\be)_{n+1} \over (n+1)! (\ga)_n }.
$$
\label{lemmgosper}
\end{lemm}

The sum (\ref{resb}) then follows from this lemma by putting
$\al=1-k-2m_1$, $\be=1+2m_2$ and $\ga=2-k+2m_2$. This completes the
proof of Proposition~\ref{propH}. 

\begin{prop}
The action of $\De(Z_+)$ on the auxiliary basis vectors is
given by
\beq
\De(Z_+) w^{j_1,j_2}_{m_1,m_2} = w^{j_1,j_2}_{m_1+1,m_2} +
w^{j_1,j_2}_{m_1,m_2+1}, 
\label{DelZw}
\eeq
where $w^{j_1,j_2}_{m_1,m_2}$ is interpreted as 0 if one of the indices
$m_i>j_i$.  
\label{propZp}
\end{prop}

This is proved by direct computation, and does not involve any
combinatorial identities. The left hand side of (\ref{DelZw}) yields
\bea
&&(1\otimes Z_+ + Z_+\otimes 1)\left( \sum_{k,l\geq 0} \sum_{n\geq 0} 
a^{m_1,m_2}_{k,l} (-h^2/4)^{n} v^{j_1}_{m_1+k+n}\otimes
v^{j_2}_{m_2+l+n} \right) \nn\\
&&=(1\otimes Z_+ + Z_+\otimes 1)\left\{ \sum_{k,l\geq 0} \left(\sum_{n\geq 0} 
a^{m_1,m_2}_{k-n,l-n} (-h^2/4)^{n}\right) v^{j_1}_{m_1+k}\otimes
v^{j_2}_{m_2+l} \right\}.
\label{tempZw}
\eea
{}From the definition (\ref{defa}) of the coefficients
$a^{m_1,m_2}_{k,l}$ it follows that 
\bea
&& \sum_{n\geq 0} a^{m_1,m_2}_{k-n,l-n} (-h^2/4)^{n} = \nn\\
&&(-1)^k (h/2)^{k+l} \sum_{n\geq 0} (b^{m_1,m_2}_{k-n,l-n} -
 b^{m_1,m_2}_{k-n-1,l-n-1} ) =
(-1)^k (h/2)^{k+l} b^{m_1,m_2}_{k,l}.
\eea
Putting this back in (\ref{tempZw}) gives
\bea
&&\sum_{k,l\geq 0} (-1)^k(h/2)^{k+l} b^{m_1,m_2}_{k,l}\; v^{j_1}_{m_1+k}
\otimes v^{j_2}_{m_2+l+1} + \nn\\
&&\qquad \sum_{k,l\geq 0} (-1)^k(h/2)^{k+l} 
b^{m_1,m_2}_{k,l} \; v^{j_1}_{m_1+k+1} \otimes v^{j_2}_{m_2+l} \nn\\
&&=\sum_{k,l\geq 0} (-1)^k(h/2)^{k+l-1} (b^{m_1,m_2}_{k,l-1} -
b^{m_1,m_2}_{k-1,l}) \; v^{j_1}_{m_1+k} \otimes v^{j_2}_{m_2+l}.
\eea
On the other hand, the rhs of (\ref{DelZw}) leads to
\bea
&&\sum_{k,l\geq 0} (a^{m_1,m_2+1}_{k,l-1} - a^{m_1+1,m_2}_{k-1,l})\;
v^{j_1}_{m_1+k} \otimes v^{j_2}_{m_2+l} = 
\sum_{k,l\geq 0} (-1)^k (h/2)^{k+l-1} \bigl(b^{m_1,m_2+1}_{k,l-1} -
\nn\\ 
&& b^{m_1,m_2+1}_{k-1,l-2}- b^{m_1+1,m_2}_{k-1,l}+ b^{m_1+1,m_2}_{k-2,l-1}
\bigr) v^{j_1}_{m_1+k} \otimes v^{j_2}_{m_2+l}.
\eea
So, it remains to show that
\beq
b^{m_1,m_2}_{k,l-1} - b^{m_1,m_2}_{k-1,l} = b^{m_1,m_2+1}_{k,l-1} - 
b^{m_1,m_2+1}_{k-1,l-2}- b^{m_1+1,m_2}_{k-1,l}+
b^{m_1+1,m_2}_{k-2,l-1}, 
\eeq
and this follows trivially from the definition (\ref{defb}) of the
$b$-coefficients. 

\begin{prop}
The action of $\De(Z_-)$ on the auxiliary basis vectors is
given by
\bea
\De(Z_-) w^{j_1,j_2}_{m_1,m_2} &=& (j_1+m_1)(j_1-m_1+1)\; 
w^{j_1,j_2}_{m_1-1,m_2} + \nn\\
&& \qquad (j_2+m_2)(j_2-m_2+1)\; w^{j_1,j_2}_{m_1,m_2-1}.
\label{DelZ-w}
\eea
\label{propZm}
\end{prop}

The proof of this property is rather long and technical. Let $h$,
$z_\pm$ be the standard basis of $sl(2)$, with
\beq
[h,z_\pm]=\pm 2 z_\pm,\qquad [z_+,z_-]=h;
\eeq
$sl(2)$ has an action on the basis elements $v^j_m$, given by the same
expressions as in (\ref{sl2vaction}). The standard comultiplication for
$sl(2)$, 
$\de(x)=x \otimes 1 + 1 \otimes x$ for every $x\in sl(2)$, induces an
action on elements of $V^{(j_1)}\otimes V^{(j_2)}$~:
\bea
&&\de(h) v^{j_1}_{m_1}\otimes v^{j_2}_{m_2} = 2(m_1+m_2)\;
v^{j_1}_{m_1}\otimes v^{j_2}_{m_2} ,\nn\\
&&\de(z_+) v^{j_1}_{m_1}\otimes v^{j_2}_{m_2} = 
v^{j_1}_{m_1+1}\otimes v^{j_2}_{m_2} + v^{j_1}_{m_1}\otimes
v^{j_2}_{m_2+1}, \\
&&\de(z_-) v^{j_1}_{m_1}\otimes v^{j_2}_{m_2} = (j_1+m_1)(j_1-m_1+1)\;
v^{j_1}_{m_1+1}\otimes v^{j_2}_{m_2} + \nn\\
&&\qquad\qquad\qquad (j_2+m_2)(j_2-m_2+1)\; v^{j_1}_{m_1}\otimes
v^{j_2}_{m_2+1}. \nn 
\eea
Thus the action of $\de(h)$ and $\de(z_+)$ on the basis
$v^{j_1}_{m_1}\otimes v^{j_2}_{m_2}$ is the same as the action of
$\De(H)$ and $\De(Z_+)$ on the auxiliary basis $w^{j_1,j_2}_{m_1,m_2}$,
by Propositions~\ref{propH} and~\ref{propZp}. We shall show that this holds 
for $Z_-$ too.
First of all, for $sl(2)$ we know how the tensor product decomposes, so
let 
\beq
v^{(j_1j_2)j}_m = \sum_{m_1+m_2=m} c^{j_1,j_2,j}_{m_1,m_2,m}\;
v^{j_1}_{m_1}\otimes v^{j_2}_{m_2} ,
\label{vvv}
\eeq
where $j\in J=\{j_1+j_2,j_1+j_2-1,\ldots, |j_1-j_2|\}$, and
$c^{j_1,j_2,j}_{m_1,m_2,m}$ are the Clebsch-Gordan coefficients in the
$v$-basis, related to the usual $su(2)$ Clebsch-Gordan coefficients by
\beq
c^{j_1,j_2,j}_{m_1,m_2,m}= {\al_{j,m}\over \al_{j_1,m_1}\al_{j_2,m_2}}
 C^{j_1,j_2,j}_{m_1,m_2,m},
\eeq
see (\ref{ve}). Then there holds
\bea
&&\de(h) v^{(j_1j_2)j}_m  = 2m \;v^{(j_1j_2)j}_m ,\nn\\
&&\de(z_+) v^{(j_1j_2)j}_m = v^{(j_1j_2)j}_{m +1}, \\
&&\de(z_-) v^{(j_1j_2)j}_m = (j+m)(j-m+1) \; v^{(j_1j_2)j}_{m-1}. \nn
\eea 
We also define
\beq
w^{(j_1j_2)j}_m = \sum_{m_1+m_2=m} c^{j_1,j_2,j}_{m_1,m_2,m} \;
w^{j_1,j_2}_{m_1,m_2},
\label{ww}
\eeq
and by the remark that 
the action of $\de(h)$ and $\de(z_+)$ on the basis
$v^{j_1}_{m_1}\otimes v^{j_2}_{m_2}$ is the same as the action of
$\De(H)$ and $\De(Z_+)$ on the auxiliary basis $w^{j_1,j_2}_{m_1,m_2}$,
it follows that
\bea
&&\de(H) w^{(j_1j_2)j}_m  = 2m \; w^{(j_1j_2)j}_m ,\nn\\
&&\de(Z_+) w^{(j_1j_2)j}_m = w^{(j_1j_2)j}_{m +1}. 
\eea
It remains to find the action of $\De(Z_-)$ on $w^{(j_1j_2)j}_m$;
let us denote it by
\beq
\De(Z_-) w^{(j_1j_2)j}_m = \sum_{j'} \sum_{m'} \mu^{j,j'}_{m,m'} \;
w^{(j_1j_2)j'}_{m'},
\eeq
with $\mu^{j,j'}_{m,m'}$ the coefficients to be determined. By acting
with the relation
\beq
\De(H)\De(Z_-)-\De(Z_-)\De(H)= -2\De(Z_-)
\eeq
on $w^{(j_1j_2)j}_m$, it follows that the coefficients
$\mu^{j,j'}_{m,m'}$ are zero unless $m'=m-1$. Write $\nu^{j,j'}_m$ for
$\mu^{j,j'}_{m,m-1}$; then we have so far
\beq
\De(Z_-) w^{(j_1j_2)j}_m = \sum_{j'} \nu^{j,j'}_{m} \;
w^{(j_1j_2)j'}_{m-1}, 
\eeq
where $j'\in J$ such that $j'\geq |m-1|$. Next, we shall use the
relation 
\beq
\De(Z_+)\De(Z_-) - \De(Z_-)\De(Z_+) = \De(H).
\label{rel}
\eeq
Acting with (\ref{rel}) on $w^{(j_1j_2)j}_m$ yields~:
\bea
\hbox{for } j\ne j'&:& \nu^{j,j'}_m = \nu^{j,j'}_{m+1}, \label{nueq1}\\
\hbox{for }j'=j & : & \nu^{j,j}_m - \nu^{j,j}_{m+1} = 2m. \label{nueq2}
\eea
In particular, by acting with (\ref{rel}) on $w^{(j_1j_2)j}_j$, one
finds that $\nu^{j,j}_j=2j$, and that $\nu^{j,j'}_j=0$ for $j'>j$. 
Now (\ref{nueq2}) is a difference equation in $m$ with boundary
condition $\nu^{j,j}_j=2j$, so it has a unique solution given by 
\beq
\nu^{j,j}_m = (j+m)(j-m+1).
\eeq
{}From (\ref{nueq1}) and $\nu^{j,j'}_j=0$ for $j'>j$ it follows that
$\nu^{j,j'}_m=0$ for all $j'>j$. Using this, and acting with
(\ref{rel}) on $w^{(j_1j_2)j}_{j-1}$ implies that
$\nu^{j,j-1}_{j-1}=0$, thus by (\ref{nueq1}) that $\nu^{j,j-1}_{m}=0$
for all $m$. Similarly, acting with (\ref{rel}) on
$w^{(j_1j_2)j}_{j-2}$ implies that $\nu^{j,j-2}_{j-2}=0$, thus by
(\ref{nueq1}) that $\nu^{j,j-2}_{m}=0$ for all $m$. One can continue
and thus show by induction that $\nu^{j,j'}_m=0$ also for all $j'<j$.
The final result is that 
\beq
\De(Z_-) w^{(j_1,j_2)j}_m = (j+m)(j-m+1)\; w^{(j_1,j_2)j}_{m-1}.
\eeq
In other words, $\De(Z_-)$ has on the basis $w^{(j_1,j_2)j}_m$ the same
action as $\de(z_-)$ on the basis $v^{(j_1j_2)j}_m$. By relations
(\ref{vvv}) and (\ref{ww}), it follows that the action of $\De(Z_-)$ on the
basis $w^{j_1,j_2}_{m_1,m_2}$ is the same as the action of $\de(z_-)$
on the basis $v^{j_1}_{m_1} \otimes v^{j_2}_{m_2}$. This proves
Proposition~\ref{propZm}.

\section{Clebsch-Gordan coefficients for $\U$}
\setcounter{equation}{0}

Let us first normalise the auxiliary basis $w^{j_1,j_2}_{m_1,m_2}$ as
follows~: 
\beq
e^{j_1,j_2}_{m_1,m_2}= w^{j_1,j_2}_{m_1,m_2}
/(\al_{j_1,m_1}\al_{j_2,m_2}) .
\eeq
Then (\ref{defw}) and (\ref{ve}) imply that
\beq
e^{j_1,j_2}_{m_1,m_2}= \sum_{k,l\geq 0} A^{m_1,m_2}_{k,l}
e^{j_1}_{m_1+k} \otimes e^{j_2}_{m_2+l},
\label{eee}
\eeq
where
\beq
A^{m_1,m_2}_{k,l} = a^{m_1,m_2}_{k,l} {\al_{j_1,m_1+k}\al_{j_2,m_2+l} \over
\al_{j_1,m_1}\al_{j_2,m_2}};
\label{defA}
\eeq
note that the $A$-coefficients depend implicitly also on $j_1$ and
$j_2$~: $A^{m_1,m_2}_{k,l}=0$ unless $m_1$ and $m_1+k$ belong to
$\{-j_1, -j_1+1,\ldots ,j_1\}$ and $m_2$ and $m_2+l$ belong to
$\{-j_2, -j_2+1,\ldots ,j_2\}$.
{}From the previous section it follows that the following relations
hold~: 
\bea
\De(H) e^{j_1,j_2}_{m_1,m_2} &=& 2(m_1+m_2)\; e^{j_1,j_2}_{m_1,m_2},\nn\\
\De(Z_+) e^{j_1,j_2}_{m_1,m_2}& &= \sqrt{(j_1-m_1)(j_1+m_1+1)}\;
e^{j_1,j_2}_{m_1+1,m_2} + \nn\\
&& \qquad \sqrt{(j_2-m_2)(j_2+m_2+1)}\; e^{j_1,j_2}_{m_1,m_2+1} ,\\
\De(Z_-) e^{j_1,j_2}_{m_1,m_2} &=& \sqrt{(j_1+m_1)(j_1-m_1+1)}\;
e^{j_1,j_2}_{m_1-1,m_2} + \nn\\
&&\qquad \sqrt{(j_2+m_2)(j_2-m_2+1)}\; e^{j_1,j_2}_{m_1,m_2-1} . \nn
\eea
Thus the action of $\De(H)$, $\De(Z_\pm)$ on $e^{j_1,j_2}_{m_1,m_2}$ is
the same as the action of $\de(h)$, $\de(z_\pm)$ on $e^{j_1}_{m_1}
\otimes e^{j_2}_{m_2}$. Consequently, we can write
\beq
e^{(j_1j_2)j}_m = \sum_{m_1+m_2=m} C^{j_1,j_2,j}_{m_1,m_2,m}\;
e^{j_1,j_2}_{m_1,m_2}, 
\label{defe}
\eeq
with $C^{j_1,j_2,j}_{m_1,m_2,m}$ the classical $su(2)$ Clebsch-Gordan
coefficients given by~\cite{Edmonds}
\bea
&&C^{j_1,j_2,j}_{m_1,m_2,m}=\de_{m_1+m_2,m} \left( (j_1+j_2-j)!
(j_1-j_2+j)! (-j_1+j_2+j)! \over (j_1+j_2+j+1)! \right)^{1/2} \nn\\
&&\times\bigl((j_1+m_1)!(j_1-m_1)!(j_2+m_2)!(j_2-m_2)!(j+m)!(j-m)!(2j+1)
\bigr)^{1/2} \nn\\
&&\times \sum_k (-1)^k/\bigl(k!(j_1+j_2-j-k)!(j_1-m_1-k)!(j_2+m_2-k)!\nn\\
&&\qquad\qquad \times (j-j_2+m_1+k)! (j-j_1-m_2+k)!\bigr).
\eea
The action of $\De(H)$ and $\De(Z_\pm)$ is then indeed
\bea
&&\De(H) e^{(j_1j_2)j}_m = 2m \;e^{(j_1j_2)j}_m ,\nn\\
&&\De(Z_\pm) e^{(j_1j_2)j}_m = \sqrt{(j\mp m)(j\pm m+1)}\;
e^{(j_1j_2)j}_{m\pm 1},
\eea
i.e.\ the same as the standard action (\ref{sl2eaction}) of $H$ and
$Z_\pm$ on a 
basis $e^j_m$. Consequently, also for $\De(X)$ and $\De(Y)$ the action
on $e^{(j_1j_2)j}_m$ is the same as the standard action of $X$ and $Y$
on a basis $e^j_m$, and the $e^{(j_1j_2)j}_m$ are genuinly ``coupled
states'' for $\U$. The decomposition of the tensor product for $\U$ is
the same as in $su(2)$~:
\beq
V^{(j_1)} \otimes V^{(j_2)} = \bigoplus_{j=|j_1-j_2|}^{j_1+j_2}
V^{(j)}. 
\eeq

{}From (\ref{eee}) and (\ref{defe}), it follows that we can write
\bea
e^{(j_1j_2)j}_m &=& \sum_{m_1+m_2=m} C^{j_1,j_2,j}_{m_1,m_2,m}
\sum_{k,l\geq 0} A^{m_1,m_2}_{k,l} \; e^{j_1}_{m_1+k} \otimes
e^{j_2}_{m_2+l} \nn\\
&=& \sum_{n_1,n_2} \Bigl( \sum_{m_1+m_2=m} C^{j_1,j_2,j}_{m_1,m_2,m} 
A^{m_1,m_2}_{n_1-m_1,n_2-m_2} \Bigr) e^{j_1}_{n_1} \otimes
e^{j_2}_{n_2} .
\eea

Thus we have
\begin{theo}
The Clebsch-Gordan coefficients for $\U$, in
\beq
e^{(j_1j_2)j}_m = \sum_{n_1,n_2} {\cal C}^{j_1,j_2,j}_{n_1,n_2,m}(h) \;
e^{j_1}_{n_1} \otimes e^{j_2}_{n_2},
\eeq
are given by
\beq
{\cal C}^{j_1,j_2,j}_{n_1,n_2,m}(h) = \sum_{m_1+m_2=m}
C^{j_1,j_2,j}_{m_1,m_2,m} A^{m_1,m_2}_{n_1-m_1,n_2-m_2},
\label{cgc}
\eeq
with $C^{j_1,j_2,j}_{m_1,m_2,m}$ the usual $sl(2)$ Clebsch-Gordan
coefficients, and $A^{m_1,m_2}_{n_1-m_1,n_2-m_2}$ determined by
(\ref{defA}) and (\ref{defa}).
\label{theoCGC}
\end{theo}

One question that naturally arises is whether these Clebsch-Gordan
coefficients satisfy an orthogonality relation. The answer is negative.
In order to have an orthogonality relation, one needs a $*$-Hopf
algebra~\cite[\S 4.1.F]{Chari}.
For $\U$ the obvious choice of the $*$-operation would be that induced
{}from $H^*=H$, $Z_\pm^*=Z_\mp$ (making the $e$-basis of $V^{(j)}$ an
orthonormal basis). However, this $*$-operation is not compatible with
the coalgebra structure.

We end this section by noting that $\U$, and its Clebsch-Gordan
coefficients, also have another interpretation. Due to the invertible
nonlinear map described in Section~2, one can identify the algebra part
of $\U$ with $U(sl(2))$. In other words, the Hopf algebra one is
dealing with has $U(sl(2))$ as algebra structure, with generators
$H$, $Z_\pm$, and defining relations~(\ref{defrelsl2}). There is no
deformation in the algebra part. The deformation comes in at the
level of the coalgebra part~: see e.g.\ (\ref{DelH}) and (\ref{DelZ}).
So, roughly speaking we are dealing with the ordinary $sl(2)$ algebra
and its usual finite dimensional irreducible representations, but
equipped with a deformed comultiplication. The deformation of the
Clebsch-Gordan coefficients then stems from this deformed coproduct.

\section{Examples and conclusion}
\setcounter{equation}{0}

The formula (\ref{cgc}) allows one to calculate the $\U$ Clebsch-Gordan
coefficients for arbitrary parameters, since the usual $sl(2)$
coefficients $C^{j_1,j_2,j}_{m_1,m_2,m}$ are known, and the
coefficients $A^{m_1,m_2}_{n_1-m_1,n_2-m_2}$ are determined in this
paper. For example, one finds
\beas
{\cal C}^{2,2,3}_{2,0,2}(h)&=&C^{2,2,3}_{2,0,2} = 1/\sqrt{2};\\
{\cal C}^{2,2,3}_{2,0,3}(h)&=& 0 ;\\
{\cal C}^{2,2,3}_{2,0,-1}(h)&=& -6h^3 C^{2,2,3}_{0,-1,-1} -4\sqrt{6}
h^3 C^{2,2,3}_{1,-2,-1} = -18 h^3/\sqrt{5}.
\eeas
More generally, suppose $m=n_1+n_2+p$, then (\ref{cgc}) becomes
\beq
{\cal C}^{j_1,j_2,j}_{n_1,n_2,m}(h) = \sum_{m_1}
C^{j_1,j_2,j}_{m_1,m-m_1,m} A^{m_1,m-m_1}_{n_1-m_1,-n_1+m_1-p}.
\label{cgc0}
\eeq
First, let $p=0$. Since $A^{m_1,m_2}_{k,l}$ is zero if $k$ or $l$ are
negative, it follows that the only contribution in (\ref{cgc0}) is for
$m_1=n_1$, and with $A^{m_1,m_2}_{0,0}=1$, it follows that ${\cal
C}^{j_1,j_2,j}_{n_1,n_2,n_1+n_2}(h)= C^{j_1,j_2,j}_{n_1,n_2,n_1+n_2}$. 
Next, suppose that $p>0$, then one can see that at least one of the
indices of $A$ in (\ref{cgc0}) is negative, thus in this case 
${\cal C}^{j_1,j_2,j}_{n_1,n_2,n_1+n_2+p}(h)=0$. Finally, suppose that
$p<0$. Now there can be a number of contributions in (\ref{cgc0}), and
by (\ref{defA}) and (\ref{defa}) they all have the same power of $h$,
namely $h^{-p}$. Thus we have the following property~:

\begin{prop}
The Clebsch-Gordan coefficients for $\U$ satisfy 
\begin{itemize}
\item if $m=n_1+n_2$ then ${\cal
C}^{j_1,j_2,j}_{n_1,n_2,m}(h) = C^{j_1,j_2,j}_{n_1,n_2,m}$ ;
\item if $m>n_1+n_2$ then ${\cal C}^{j_1,j_2,j}_{n_1,n_2,m}(h) =0$;
\item if $m<n_1+n_2$ then ${\cal C}^{j_1,j_2,j}_{n_1,n_2,m}(h)$ is a
monomial in $h^{n_1+n_2-m}$.
\end{itemize}
\end{prop}

This proposition implies that the Clebsch-Gordan coefficients ${\cal
C}^{j_1,j_2,j}_{n_1,n_2,m}(h)$ are simple deformations of the $sl(2)$
Clebsch-Gordan coefficients in the sense that for $h=0$ one has ${\cal
C}^{j_1,j_2,j}_{n_1,n_2,m}(0) = C^{j_1,j_2,j}_{n_1,n_2,m}$. 
As is well known, $sl(2)$ or $su(2)$ representations play a fundamental
role in many physical models and theories, and so do their Clebsch-Gordan
coefficients. 
It would be interesting to investigate whether the present
Clebsch-Gordan coefficients of $\U$ still have a physical
interpretation in corresponding deformed models or theories.
The explicit formula given here, together with the properties
mentioned, should prove to be very helpful in such an investigation.

\section*{Appendix}
\renewcommand{\theequation}{A.\arabic{equation}}
\setcounter{equation}{0}

Here, we give proofs of the combinatorial identities in 
Lemma~\ref{lemmYcoef}, \ref{lemmzeilb} and~\ref{lemmgosper}.
The terms appearing in the sums of these lemmas are
hypergeometric terms~\cite{AB}. The sums in Lemma~\ref{lemmYcoef} 
and~\ref{lemmzeilb}  are definite sums,
i.e.\ the summation limit also appears in the summand. For such sums,
Zeilberger's algorithm~\cite[Chapter~6]{AB}, also known as the method of
creative telescoping~\cite{Zeilberger} can be used to find a
recurrence relation. A Mathematica implementation {\tt Zb} of
Zeilberger's algorithm can for example be found in the package {\tt
Zb.m}~\cite{Paule}. Considering the sum of Lemma~\ref{lemmYcoef},
\beq
f(k) = \sum_{n=1}^{k} (-1)^n {(1/2)_n \over n!} \left({k-1 \atop
n-1}\right), 
\eeq
Zeilberger's algorithm yields the following recurrence relation for
$f(k)$~: 
\beq
2(k+1)f(k+1)=(2k-1) f(k).
\eeq
with the initial condition $f(1)=-1/2$. The closed form expression for
$f(k)$ then easily follows and is given in Lemma~\ref{lemmYcoef}.

Lemma~\ref{lemmzeilb} is similar, with
\beq
f_n(s) = \sum_{k=1}^{s-1} {(2k-2)!\over k!(k-1)!}
{(2s-2k-2)!\over(s-k)!(s-k-1)!} k^n,
\eeq
where $n=0$, 1 or 2. The recurrence relations obtained by Zeilberger's
algorithm (or by the program {\tt Zb}) read
\bea
&&4(s-1)f_0(s)-(s+1) f_0(s+1)+{2(2s-2)!\over s!(s-1)!}=0,\\
&&4(s-1)f_1(s)-s f_1(s+1)+{(2s-2)!\over (s-1)!^2}=0,\\
&&4f_2(s)- f_2(s+1)+{(2s-2)!\over s! (s-2)!}=0.
\eea
With the initial conditions $f_n(2)=1$, the closed form expressions for
$f_n(s)$ given in Lemma~\ref{lemmzeilb} are deduced from these recurrence 
relations.

The statement in Lemma~\ref{lemmgosper} is rather different, in the sense 
that 
this is an indefinite summation (i.e.\ the upper limit does not appear in
$g(l)$). The term $g(l)$ is again a hypergeometric term however. For
such summations Gosper's algorithm~\cite{AB,Gosper} decides whether the
sum can be 
written in closed form, and also gives the closed form if it exists.
The summation formula in Lemma~\ref{lemmgosper} is a direct output of the 
Mathematica program {\tt Gosper} of the package {\tt Zb.m}. 


\end{document}